# FIRST ORDER DIFFERENTIAL EQUATIONS INDUCED BY THE INFINITE SERVERS QUEUE WITH POISSON ARRIVALS TRANSIENT BEHAVIOR PROBABILITY DISTRIBUTION PARAMETERS STUDY AS TIME FUNCTIONS

FERREIRA Manuel Alberto M. (PT)

**Abstract.** The M|G|∞ queue system state transient probabilities, considering the time origin at the beginning of a busy period, mean and variance monotony as time functions is studied. These studies, for which results it is determinant the hazard rate function service time length, induce the consideration of two differential equations, one related with the mean monotony study and another with the variance monotony study, which solutions lead to some particular service time distributions, for which those parameters present specific behaviors as time functions.

## 1   Introduction

The M|G|∞ queue customers arrive at this queue system according to a Poisson process at rate $\lambda$ and receive a service which time length is a positive random variable with distribution function $G(.)$ and mean $\alpha$. Upon their arrival, they find immediately an available server. Each customer service time length is independent from the other customers' services time length and from the arrivals process. The traffic intensity is a no dimensional parameter given by $\rho = \lambda\alpha$.

Call $N(t)$ the number of occupied servers (or the number of customers being served) at instant $t$, in a M|G|∞ queue system.

Be $p_{1'n} = P[N(t) = n | N(0) = 1']$, $n = 0,1,2,...$, meaning $N(0) = 1'$ that the time origin is an instant at which a customer arrives at the system, jumping the number of customers from 0 to 1. That is: an instant at which a busy period begins. Then, see (Ferreira and Andrade, 2009),

$$p_{1'0}(t) = p_{00}(t)G(t)$$
$$p_{1'n}(t) = p_{0n}(t)G(t) + p_{0n-1}(t)(1-G(t)), \ n = 1,2,... \tag{1.1}$$

where

$$p_{0n}(t) = \frac{\left(\lambda \int_0^t [1-G(v)]dv\right)^n}{n!} e^{-\lambda \int_0^t [1-G(v)]dv}, \quad n = 0,1,2,\ldots, \tag{1.2}$$

being here the time origin an instant at which the queue system is customers empty, see again (Ferreira and Andrade, 2009).

It is easy to check that

$$\lim_{t \to \infty} p_{1'n}(t) = \frac{\rho^n}{n!} e^{-\rho}, \quad n = 0,1,2,\ldots \tag{1.3}$$

such as $\lim_{t \to \infty} p_{0n}(t) = \frac{\rho^n}{n!} e^{-\rho}, \quad n = 0,1,2,\ldots$

Denoting $\mu(1',t)$ the distributions given by (1.1) mean value:

$$\mu(1',t) = 1 - G(t) + \lambda \int_0^t [1 - G(v)]dv \tag{1.4}$$

And denoting $V(1',t)$ the distribution defined by (1.1) variance, it is obtained

$$V(1',t) = \lambda \int_0^t [1 - G(v)]dv + G(t) - G^2(t). \tag{1.5}$$

The main target is to study $\mu(1',t)$ and $V(1',t)$ monotony as time functions. Then two first order differential equations, induced by this study, are considered allowing very interesting conclusions. In this work are continued and deepened the works of Ferreira (2005 and 2014).

## 2  Study of $\mu(1',t)$ as time function

Call $h(t)$ the service time length hazard rate function. It is given by, see Ross (1983),

$$h(t) = \frac{g(t)}{1-G(t)} \tag{2.1}$$

where $g(t) = \frac{dG(t)}{dt}$. It is interpreted as the rate at which the services end.

**Proposition 2.1**

If $G(t) < 1, t > 0$ continuous and differentiable and

$$h(t) \leq \lambda, \ t > 0 \tag{2.2}$$

$\mu(1',t)$ is non-decreasing.

**Dem.:**

It is enough to note, according to (1.4), that $\frac{d}{dt}\mu(1',t) = (1-G(t))(\lambda - h(t))$. ∎

**Obs.:**

-The expression (2.2) meaning is that if the rate at which the services end is lesser or equal than the customers arrivals rate, $\mu(1',t)$ is non-decreasing.

-For the M|M|∞ system, that is M|G|∞ queue with exponentially distributed service times length, (2.2) is equivalent to

$$\rho \geq 1 \qquad (2.3)$$

- $\lim_{t\to\infty} \mu(1',t) = \rho$, as expected. ∎

Defining $\beta(.)$ as

$$\beta(t) = h(t) - \lambda \qquad (2.4)$$

The first order differential equation (2.4) solution is the following collection of service time distributions:

$$G(t) = 1 - (1-G(0))e^{-\lambda t - \int_0^t \beta(u)du}, t \geq 0, \frac{\int_0^t \beta(u)du}{t} \geq -\lambda \qquad (2.5)$$

So, evidently,

**Proposition 2.2**

If $\beta(t) = 0$

$$G(t) = 1 - (1-G(0))e^{-\lambda t}, \quad t \geq 0 \qquad (2.6)$$

and $\mu(1',t) = 1 - G(0) = \rho, \quad t \geq 0$. ∎

**Proposition 2.3**

If $\beta(t) = -\lambda$

$$G(t) = 1, t \geq 0 \qquad (2.7)$$

and $\mu(1',t) = 0, \quad t \geq 0$. ∎

For some service time length distributions:

-Deterministic with value $\alpha$ (**M|D|∞** system)

$$\mu(1',t)=\begin{cases}1+\lambda t, & t<\alpha \\ \rho, & t\geq\alpha\end{cases} \quad (2.8)$$

-Exponential (**M|M|∞** system)

$$\mu(1',t)=\rho+(1-\rho)e^{-\frac{t}{\alpha}} \quad (2.9)$$

$$-G(t)=1-\frac{(1-e^{-\rho})(\lambda+\beta)}{\lambda e^{-\rho}(e^{(\lambda+\beta)t}-1)+\lambda}, t\geq 0, -\lambda\leq\beta\leq\frac{\lambda}{e^{\rho}-1}$$

$$\mu(1't)=\frac{(1-e^{-\rho})(\lambda+\beta)}{\lambda e^{-\rho}(e^{(\lambda+\beta)t}-1)+\lambda}+\rho-\log\left(1+(e^{\rho}-1)e^{-(\lambda+\beta)t}\right) \quad (2.10)$$

For this collection of service time distributions, the busy period time length is exponentially distributed with an atom at the origin, see (Ferreira and Andrade, 2009) and Ferreira (2016 and 2016a)**:**

$$B^{\beta}(t)=1-\frac{\lambda+\beta}{\lambda}(1-e^{-\rho})e^{-e^{-\rho}(\lambda+\beta)t}, t\geq 0, -\lambda\leq\beta\leq\frac{\lambda}{e^{\rho}-1} \quad (2.11)$$

**Note:**

-For a M/G/∞ queue, if the service time length distribution function belongs to the collection, see again (Ferreira and Andrade, 2009) and Ferreira (2016 and 2016a),

$$G(t)=1-\frac{1}{\lambda}\frac{(1-e^{-\rho})e^{-\lambda t-\int_0^t\beta(u)du}}{\int_0^{\infty}e^{-\lambda w-\int_0^w\beta(u)du}dw-(1-e^{-\rho})\int_0^t e^{-\lambda w-\int_0^w\beta(u)du}dw},$$

$$t\geq 0, -\lambda\leq\frac{\int_0^t\beta(u)du}{t}\leq\frac{\lambda}{e^{\rho}-1},$$

which is solution for the Riccati equation

$$\frac{dG(t)}{dt} = -\lambda G^2(t) - (\beta(t) - \lambda)G(t) + \beta(t),$$

the busy period time length distribution function is

$$B(t) = \left(1 - (1-G(0))\left(e^{-\lambda t - \int_0^t \beta(u)du} + \lambda \int_0^t e^{-\lambda w - \int_0^w \beta(u)du}\,dw\right)\right) *$$

$$* \sum_{n=0}^{\infty} \lambda^n (1-G(0))^n \left(e^{-\lambda t - \int_0^t \beta(u)du}\right)^{*n}, -\lambda \leq \frac{\int_0^t \beta(u)du}{t} \leq \frac{\lambda}{e^\rho - 1}.$$

If $\beta(t) = \beta$ (constant) it results (2.11). ∎

## 3 Study of $V(1',t)$ as time function

**Proposition 3.1**

If $\frac{1}{2} < G(t) < 1, t > 0$ continuous and differentiable and

$$h(t) \leq \frac{\lambda}{2G(t)-1} \tag{3.1}$$

$V(1',t)$ is non-decreasing.

**Dem.:**
According to (1.5), $\frac{d}{dt} V(1',t) = \lambda(1 - G(t)) + g(t) - 2G(t)g(t) = \lambda(1 - G(t)) + g(t)(1 - 2G(t)) = (1 - G(t))(h(t)(1 - 2G(t)) + \lambda)$. So that $\frac{d}{dt} V(1',t) \geq 0$ it is enough to have $(h(t)(1 - 2G(t)) + \lambda) \geq 0 \Leftrightarrow h(t)(2G(t) - 1) - \lambda \leq 0 \Leftrightarrow h(t) \leq \frac{\lambda}{2G(t)-1}$ since $G(t) > \frac{1}{2}$. ∎

**Obs.:**

- $\lim_{t \to \infty} V(1',t) = \rho$, as expected. ∎

Define $\beta(.)$ as $\beta(t) = h(t) - \frac{\lambda}{2G(t)-1}$. It results the following differential equation in G(.):

$$\frac{dG(t)}{dt} = \left(\beta(t) + \frac{\lambda}{2G(t)-1}\right)(1 - G(t)), t \geq 0 \tag{3.2}$$

Give (3.2) the following form

$$(\lambda + \beta(t)(2G(t) - 1))dt + \left(2 - \frac{1}{1-G(t)}\right)dG = 0 \tag{3.3}$$

that can be turned in

$$(\lambda + \beta(2G(t) - 1))dt + \left(2 - \frac{1}{1-G(t)}\right)dG = 0 \qquad (3.4)$$

for $\beta(t) = \beta$ (const.). So $\frac{\partial}{\partial G}(\lambda + \beta(2G - 1)) = 2\beta$ and $\frac{\partial}{\partial t}\left(2 - \frac{1}{1-G}\right) = 0$. Then it is possible to consider the following situations:

A) $\boldsymbol{\beta = 0}$

Equation (3.2) is an exact differential equation which solution is

$$\frac{1-G(t)}{1-G(0)} e^{2(G(t)-G(0))} = e^{-\lambda t}, t \geq 0 \qquad (3.5)$$

**Proposition 3.2**

If $G(.)$ is implicitly defined as

$$\frac{1-G(t)}{1-G(0)} e^{2(G(t)-G(0))} = e^{-\lambda t}, t \geq 0 \qquad (3.6)$$

$V(1',t) = \rho, t \geq 0.$ ∎

**Obs.:**

-The density associated to (3.6) is

$$g(t) = -\frac{\lambda e^{-\lambda t}(1-G(0))}{(1-2G(t))e^{2(G(t)-G(0))}} . \qquad (3.7)$$

-After (3.7), denoting $S$ the associated random variable, it is easy to see that, with $G(0) > \frac{1}{2}$,

$$\frac{(1-G(0))n!e^{-2(1-G(0))}}{\lambda^n} \leq E[S^n] \leq \frac{(1-G(0))n!}{(2G(0)-1)\lambda^n}, n = 1,2,... . \blacksquare \qquad (3.8)$$

B) $\boldsymbol{\beta \neq 0}$

Now equation (3.2) is not an exact differential equation but, using the integrating factor

$$\mu = \frac{1}{\beta(2G(t)-1)+\lambda} \tag{3.9}$$

that is a solution for the differential equation

$$\frac{d\ln\mu}{dG} = \frac{\frac{\partial}{\partial t}\left(2-\frac{1}{1-G}\right)-\frac{\partial}{\partial G}(\lambda+\beta(2G-1))}{\lambda+\beta(2G-1)} = \frac{0-2\beta}{\lambda+\beta(2G-1)} = -\frac{2\beta}{\beta(2G-1)+\lambda} \tag{3.10}$$

it becomes an exact differential equation which solution is

$$\left(\frac{1-G(t)}{1-G(0)}\right)^{\frac{1}{\beta+\lambda}} \left|\frac{2\lambda G(t)-\beta+\lambda}{2\lambda G(0)-\beta+\lambda}\right|^{\frac{\lambda}{\beta(\beta+\lambda)}} = e^{-t} . \tag{3.11}$$

It is interesting to make $\beta = \lambda$ in (3.11). In fact it becomes

$$\left(\frac{1-G(t)}{1-G(0)}\right)\left(\frac{G(t)}{G(0)}\right) = e^{-2\lambda t}$$

that lead to

$$G(t) = \frac{1+\sqrt{1-4(1-G(0))G(0)e^{-2\lambda t}}}{2} . \tag{3.12}$$

It is easy to check that in (3.12) $G(\infty) = 1$ and $G(0) = \frac{1}{2}$.

**Note:**

- Cases A) and B) may be resolved using the separation variables technic, maybe easier. The proceeding used comes from an attempt to get solutions for the general case with $\beta(t)$ non constant. But, in this situation, the second member of (3.10) does not depend only on $G$ what makes this proceeding to determine the integrating factor unviable.

- $G(t) = 1, t \geq 0$ is a trivial solution of equation (3.2). And, for this case, it is easy to check after (1.5) that $V(1',t)=0, t \geq 0$. ∎

For some particular service time length distributions:

-Deterministic with value $\alpha$ (**M|D|∞** system)

$$V(1',t)=\begin{cases} \lambda t, t < \alpha \\ \\ \rho, t \geq \alpha \end{cases} . \tag{3.13}$$

-Exponential (**M|M|∞** system)

$$V(1',t) = \rho\left(1-e^{-t/\alpha}\right) + e^{-\frac{t}{\lambda}} + e^{-\frac{2t}{\alpha}}, \tag{3.14}$$

$$-G(t) = 1 - \frac{\left(1-e^{-\rho}\right)(\lambda+\beta)}{\lambda e^{-\rho}\left(e^{(\lambda+\beta)t}-1\right)+\lambda}, t \geq 0, -\lambda \leq \beta \leq \frac{\lambda}{e^{\rho}-1}$$

$$V(1',t) = \rho - \log\left(1+\left(e^{\rho}-1\right)e^{-(\lambda+\beta)t}\right) + \frac{\left(1-e^{-\rho}\right)(\lambda+\beta)}{\lambda e^{-\rho}\left(e^{(\lambda+\beta)t}-1\right)+\lambda} +$$

$$\left(\frac{\left(1-e^{-\rho}\right)(\lambda+\beta)}{\lambda e^{-\rho}\left(e^{(\lambda+\beta)t}-1\right)+\lambda}\right)^2. \tag{3.15}$$

## 4  Conclusions

It was possible to study $\mu(1',t)$ and $V(1',t)$, as time functions, playing there an important role the service time length hazard rate function. This circumstance leads to differential equations (2.4) and (3.2). The first was completely solved and the second only in some particular cases. As consequence, some service time length distributions, for which those parameters have particular behaviors, as time functions, were determined.

**Current address**
**Manuel Alberto M. Ferreira, Professor Catedrático**
INSTITUTO UNIVERSITÁRIO DE LISBOA (ISCTE-IUL)
BRU – IUL, ISTAR-IUL

Av. das Forças Armadas, 1649-026 Lisboa, Portugal
Tel.: + 351 21 790 37 03. FAX: + 351 21 790 39 41,
E-mail: manuel.ferreira@iscte.pt